\documentclass[12pt,twoside]{article}
\usepackage{amssymb, amsmath}

 \usepackage{graphicx}

\newtheorem{Theorem}{Theorem}

\newtheorem{Corollary}{Corollary}

\begin{document}

\title{On an Extension of the Concept of Slowly Varying Function with Applications to Large Deviation
Limit Theorems}
\author{
A.\,A. Borovkov$^1$\ and K.\,A. Borovkov$^2$}

\footnotetext[1]{Sobolev Institute of Mathematics, Ac.~Koptyug pr.~4, 630090
Novosibirsk, Russian Federation (e-mail: borovkov@math.nsc.ru) and Novosibirsk State
Univerity.}

\footnotetext[2]{Department of Mathematics \& Statistics, The University of Melbourne,
Parkville 3010, Australia (e-mail: borovkov@unimelb.edu.au).}

\date{}

\maketitle

\begin{abstract}
Karamata's integral representation for slowly varying functions is extended to a broader
class of the so-called $\psi$-locally constant functions, i.e.\ functions $f(x)>0$
having the property that, for a given non-decreasing function~$\psi (x)$ and any fixed
$v$, ${f (x+v\psi(x) )}/{f(x)}\to 1 $ as $x\to\infty$. We   consider applications of
such functions to  extending known results on large deviations of sums of random
variables with regularly varying distribution tails.

\medskip
{\it Key words and phrases:} slowly varying function; locally constant function; large
deviation probabilities; random walk.

\medskip{\em AMS Subject Classifications 2000:} 26A12, 60F10.
\end{abstract}

\section{Introduction}

\medskip
Let  $L(x)$ be a slowly varying function  (s.v.f.), i.e.\ a positive measurable function
such that, for any fixed   $v\in(0,\infty)$ holds $L(v x)\sim L(x)$ as $x\to\infty$:
\begin{equation}
 \label{ap1}
\lim_{x\to\infty}\frac{L(v x)}{L(x)}=1.
\end{equation}
Among the most important and often used results on s.v.f.'s are the Uniform
Conver\-gence Theorem (see property~{\bf (U)} below) and the Integral Representation
Theorem (property~{\bf (I)}), the latter result essentially relying on the former. These
theorems, together with their proofs, can be found e.g.\ in the
monographs~\cite{bingham} (Theorems~1.2.1 and~1.3.1) and~\cite{bor_bor} (see~\S~1.1).

\medskip
{\bf (U)}\ {\it For any fixed $0<v_1<v_2<\infty$, convergence  \eqref{ap1} is uniform in
$v\in [v_1,v_2]$.}

\medskip
{\bf (I)} {\it A function $L(x)$ is an s.v.f.\ iff when the following representation
holds true:
\begin{equation}\label{ap2}
L(x)=c(x)\exp\bigg\{\int_1^x\frac{\varepsilon(t)}{t}\,dt\bigg\},\qquad x\geq 1,
\end{equation}
where the functions $c(t)>0$ and $\varepsilon(t)$ are measurable, $c(t)\to
c\in(0,\infty)$ and $\varepsilon(t)\to 0$ as~$t\to\infty$}.
\medskip

The concept of a s.v.f.\ is closely related to that of a regularly varying function
(r.v.f.) $R(x)$, which is specified by the relation
$$
R(x)=x^\alpha L(x),\qquad \alpha\in\mathbb{R},
$$
where $L$ is an s.v.f.\ and $\alpha$ is called the index of the r.v.f.~$R(x)$. The class
of all r.v.f.'s we will denote by~$\mathcal{R}$.

R.v.f.'s are characterised by the relation
\begin{equation}\label{ap-3}
\lim_{x\to\infty}\frac{R(vx)}{R(x)}=v^\alpha,\quad v\in(0,\infty).
\end{equation}
For them, convergence \eqref{ap-3} is also uniform in $v$ on compact intervals, while
represen\-tation~\eqref{ap2} holds for r.v.f.'s  with $\varepsilon(t)\to\alpha$ as
$t\to\infty$.

In Probability Theory there exists a whole class of limit theorems on large deviations
of sums of random variables whose distributions $\mathbf{F}$ have the property that
their right tails $F_+(x):=\mathbf{F}\big([x,\infty)\big)$ are r.v.f.'s. The following
assertion (see e.g.\ Theorem~4.4.1 in~\cite{bor_bor}) is a typical representative of
such results. Let $\xi, \xi_1, \xi_2,\ldots$ be independent identically distributed
random variables, $\mathbf{E}\xi=0$, $\mathbf{E}\xi^2<\infty$,
$S_n:=\sum\limits_{k=1}^n\xi_k$ and $\overline{S}_n:=\max\limits_{k\leq n}S_k$.

\medskip\noindent
{\bf Theorem~A}\ {\em If\, $F_+(t)=\mathbf{P}(\xi\geq t)$ is an r.v.f.\ of index
$\alpha<-2$ then, as $x\to\infty$, $ x (n \ln n)^{-1/2}\to\infty,$ one has   }
\begin{equation}\label{ap-4}
\mathbf{P}(S_n\geq x)\sim n F_+(x),\qquad \mathbf{P}(\overline{S}_n\geq x)\sim nF_+(x).
\end{equation}

\medskip Similar assertions hold true under the assumption that the distributions of the scaled sums~$S_n$
tend to a stable law as $n\to\infty$ (see Chapters~2,\,3 in~\cite{bor_bor}).

There arises the natural question of how essential the condition $F_+\in \mathcal{R}$ is
for relations~\eqref{ap-4} to hold. It turns out that this condition can be
substantially relaxed.

The aim of the present paper is to describe and study classes of functions that are
wider than~$\mathcal{R}$ and have the property that the condition that $F_+$ belongs to
such a class, together with some other natural conditions, would ensure the validity of
limit laws of the form~\eqref{ap-4}.

In Section~2 of the present note we give the definitions of the above-mentioned broad
classes of functions which we call asymptotically $\psi$-locally constant functions. The
section also contains assertions in which conditions sufficient for
relations~\eqref{ap-4} are given in terms of these functions. Section~3 presents the
main results on characterisation of asymptotically $\psi$-locally constant functions.
Section~4 contains the proofs of these results.

\section{The definitions of asymptotically locally constant functions. Applications to
limit theorems on large deviations}

Following \S~1.2 in~\cite{bor_bor}, we will call a positive function $g(x)$ an {\em
asymptotically locally constant function} (l.c.f.) if, for any
fixed~$v\in(-\infty,\infty)$,
\begin{equation}\label{ap3}
\lim_{x\to\infty}\frac{g(x+v)}{g(x)}=1
\end{equation}
(the function  $g(x)$, as all the other functions appearing in the present note, will be
assumed measurable;   assumptions of this kind will be omitted for brevity's sake).

If one puts $x:=\ln y$, $v:=\ln u$, then $g(x+v)=g(\ln (yu))$, so that the composition
$L=g\circ\ln$ will be an s.v.f.\ by virtue of~\eqref{ap3} and~\eqref{ap1}. From here and
the equality $g(x)=L(e^x)$ it follows that an l.c.f.~$g$ will have the following
properties:

\medskip

{\bf($\boldsymbol{\rm U_1}$)} {\em For any fixed $-\infty<v_1<v_2<\infty$, convergence
\eqref{ap3} is uniform in $v\in[v_1,v_2]$}.

\medskip
{\bf ($\boldsymbol{\rm I_1}$)} {\em  A function $g(x)$ is an l.c.f.\ iff it admits a
representation of the form
\begin{equation}
 \label{ap3'}
g(x)=c(x)\exp\bigg\{\int_1^{e^x}\frac{\varepsilon(t)}{t}\,dt\bigg\},\qquad x\geq 1,
\end{equation}
where $c(t)$ and $\varepsilon(t)$ have the same properties as in~{\bf(I)}}.

\medskip

Probability distributions $\mathbf{F}$ on $\mathbb{R}$ such that $F_+ (t):
=\mathbf{F}\big([t,\infty)\big)$ are  l.c.f.'s are sometimes referred to as long-tailed
distributions, or class $\mathcal{L}$ distributions. Such distribu\-tions often appear
in papers on limit theorems for sums of random variables with ``heavy tails". Examples
of l.c.f.'s  are provided by r.v.f.'s and functions of the form $\exp\{x^\alpha L(x)\}$,
where $L$ is an s.v.f., $\alpha\in(0,1)$.

It is not hard to see that, by virtue of property {\bf ($\boldsymbol{\rm U_1}$)},
definition~\eqref{ap3} of an l.c.f.\ is equivalent to the following one: for any fixed
$v\in(-\infty, \infty)$ and function $v(x)\to v$ as $x\to\infty$, one has
\begin{equation}
\label{ap4}
\lim_{x\to\infty}\frac{g (x+v(x) )}{g(x)}=1.
\end{equation}

Now we will consider a broader concept, which includes both s.v.f.'s and l.c.f.'s as
special cases.

Let  $\psi(t)>1$ be a fixed non-decreasing function.

\medskip\noindent
{\bf Definition 1} (see also Definition~1.2.7 in~\cite{bor_bor}) A function  $g(x)>0$ is
said to be  an {\em asymptotically  $\psi$-locally constant function} ($\psi$-l.c.f.)
if, for any fixed $v\in(-\infty,\infty)$ such that $x+v\psi(x)\geq cx$ for some~$c>0$
and all large enough~$x$, one has
\begin{equation}
 \label{ap5}
\lim_{x\to \infty}\frac{g (x+v\psi(x) )}{g(x)}=1.
\end{equation}

If $\psi(x)\equiv 1$ then the class of $\psi$-l.c.f.'s coincides with the class of
l.c.f.'s. If $\psi(x)\equiv x$ then the class of $\psi$-l.c.f.'s coincides with the
class of  s.v.f.'s. Now if   $\psi(x)\to\infty$ and $\psi(x)=o(x)$ as $x\to\infty$, then
the class of $\psi$-l.c.f.'s occupies, in a certain sense, an intermediate (in terms of
the zone where its functions are locally constant) place  between the classes of
s.v.f.'s and l.c.f.'s.

Clearly, all the functions from  $\mathcal{R}$ are $\psi$-l.c.f.'s for any function
$\psi(x)=o(x)$.

We will also  need the following

\medskip\noindent
{\bf Definition 2} (see also Definition~1.2.20 in~\cite{bor_bor}) \ We will call a
function   $g$ an {\em  upper-power function\/} if it is an l.c.f.\  and, for any
$p\in(0,1),$ there exists a constant $c(p)$, $\inf\limits_{p\in(p_1,1)}c(p)>0$ for
any~$p_1>0$, such that
$$
g(t)\geq c(p)g(pt).
$$
It is clear that all r.v.f.'s are upper-power functions.

\medskip

The concepts of  $\psi$-l.c.f.'s and upper-power function enable one to substantially
extend the assertion of Theorem~A. It is not hard to derive from Theorem~4.8.1
in~\cite{bor_bor} the following result.

Let  $h(v)>0$ be a non-decreasing function such that $h(v)\gg \sqrt{v\,\ln v}$ as
$v\to\infty$. Such a function always has a generalised inverse $h^{(-1)}(t):=\inf
\{v:\,h(v)\geq t \}$.

\medskip\noindent
{\bf Theorem~B} \ {\em  Let\/  $\mathbf{E}\xi=0$, $\mathbf{E}\xi^2<\infty$. Assume that
the following conditions are satisfied:

\smallskip
$1)$ $F_+(t)\leq V(t)=t^{\alpha} L(t)$, where $\alpha<-2$ and $L$ is an s.v.f.
\smallskip

$2)$ The function $F_+(t)$ is upper-power and a $\psi$-l.c.f.\ for
$\psi(t)=\sqrt{h^{(-1)}(t)}$.
\smallskip

\noindent Then relations $\eqref{ap-4}$ hold true provided that $x\to\infty$,  $x\geq
h(n)$ and
\begin{equation}
 \label{ap-5}
nV^2(x)=o\big(F_+(x)\big).
\end{equation}
}

In particular, if $x=h(n)\sim cn^\beta$ as $n\to\infty$, $\beta>1/2$, then one can put
$\psi(t):=t^{1/2\beta}$ ($\psi(t):=\sqrt{t}$ if $x\sim cn$).

Condition~\eqref{ap-5} is always met provided that $F_+(t)\geq cV(t)t^{-\varepsilon}$
for some~$\varepsilon>0$, $\varepsilon<-\alpha-2$ and $c={\rm const}$. Indeed, in this
case, for $x\geq\sqrt{n}$, $x\to\infty$,
$$
nV^2(x)\leq c^{-1}x^{2+\varepsilon}V(x)F_+(x)=o\big(F_+(x)\big).
$$

Now consider the case when $\mathbf{E}\xi^2=\infty$. Let, as before,  $V(t)=t^{\alpha}
L(t)$ be an r.v.f., $\sigma(v):=V^{(-1)}(1/v))$. Observe that $\sigma(v)$ is also an
r.v.f.\ (see e.g.\  Theorem~1.1.4 in~\cite{bor_bor}). Further, let $h(v)>0$ be a
non-decreasing function such that $h(v)\gg\sigma(v)$ as $v\to\infty$. Employing
Theorem~4.8.6 in~\cite{bor_bor} (using this opportunity, note that there are a couple of
typos in the formulation of the  theorem: the text ``with
$\psi(t)=\sigma(t)=V^{(-1)}(1/t)$" should be omitted, while the condition
``$x\gg\sigma(n)$" must be replaced with ``$x\gg\sigma(n)=V^{(-1)}(1/n)$") it is not
difficult to establish the following result.

\medskip\noindent
{\bf Theorem~C}\ {\em  Let\/ $\mathbf{E}\xi=0$ and the following conditions be met:
\smallskip

$1)$ $F_+(t)\leq V(t)=t^\alpha L(t)$, where $-\alpha\in (1,2)$ and $L$ is an s.v.f.
\smallskip

$2)$ $\mathbf{P}(\xi<-t)\leq cV(t)$ for all $t>0$.
\smallskip

$3)$ The function  $F_+$ is upper-power and a $\psi$-l.c.f.\ for
$\psi(t)=\sigma\big(h^{(-1)}(t)\big)$.
\smallskip

\noindent Then relations $\eqref{ap-4}$ hold true provided that $x\to\infty$,  $x\geq
h(n)$ and relation~\eqref{ap-4} is satisfied.}

\medskip

If, for instance, $V(t)\sim c_1 t^\alpha$ as $t\to\infty$, $x\sim c_2n^\beta$ as
$n\to\infty$, $c_i={\rm const}$, $i=1,2$, and $\beta>-1/\alpha$, then one can put
$\psi(t):=t^{-1/(\alpha\beta)}$.

Condition~\eqref{ap-5} of Theorem~C is always met provided that $x\geq
n^{\delta-(1/\alpha)}$, $F_+(t)\geq cV(t)t^{-\varepsilon}$ for some $\delta>0,$ and
$\varepsilon<  \alpha^2\delta/(1-\alpha\delta)$. Indeed, in this case $n\leq x^{ -\alpha
/(1-\alpha\delta)}$ and
$$
nV^2(x)\leq c^{-1}F_+(x)x^{\varepsilon -\alpha /(1-\alpha\delta)}V(x)=o\big(F_+(x)\big).
$$
Note also that the conditions of  Theorems~B and~C do not stipulate that $n\to\infty$.

The proofs of Theorems~B and C basically consist in verifying, for the indicated choice
of functions $\psi$, the conditions of Theorems~4.8.1 and~4.8.6 in~\cite{bor_bor},
respectively. We will omit them.

It is not hard to see (e.g.\ from the representation theorem on p.~74 in~\cite{bingham})
that Theorems~B and~C include, as special cases, situations when the right tail of
$\mathbf{F}$ satisfies the condition of {\em  extended regular variation}, i.e.\ when,
for some $0<\alpha_1\leq\alpha_2<\infty$ and any $b>1$,
\[
b^{-\alpha_2}\leq\liminf_{x\to\infty}\frac{F_+(bx)}{F_+(x)}\leq
\limsup_{x\to\infty}\frac{F_+(bx)}{F_+(x)}\leq b^{-\alpha_1}.
\]
Under the assumption that the random variable $\xi=\xi'-\mathbf{E}\xi'$ was obtained by
centering a non-negative random variable $\xi'\ge 0$, the former of the asymptotic
relations \eqref{ap-4} was established in the above-mentioned case
in~\cite{Cline-Hsing}. One could also mention here some further efforts aimed at
extending the conditions of Theorem~A that ensure the validity of~\eqref{ap-4}, see
e.g.~\cite{Cline,Ng}.

In conclusion of this section, we will make a remark showing that the presence of the
condition that $F_+(t)$ is a $\psi$-l.c.f.\ in Theorems~B and~C is quite natural and
also that is indicates that any further extension of this condition in the class of
``regular enough'' functions is hardly possible. If we turn, say, to the proof of
Theorem~4.8.1 in~\cite{bor_bor}, we will see that when $x\sim cn$, the main term in the
asymptotic representation for $\mathbf{P}(S_n\geq x)$ is given by
\[
n\int_{-N\sqrt{n}}^{N\sqrt{n}}\mathbf{P}(S_{n-1}\in dt)F_+(x-t),
\]
where $N\to\infty$ slowly enough as $n\to\infty$. It is clear that, by virtue of the
Central Limit Theorem, the integral in this expression is asymptotically equivalent to
$F_+(t)$ (implying that the former relation in~\eqref{ap-4} will hold true), provided
that $F_+(t)$ is a $\psi$-l.c.f.\ for $\psi(t)=\sqrt{t}$.

Since $\mathbf{E}S_{n-1}=0$, one might try to obtain such a result in the case when то
$F_+\big(t+v\psi(t)\big)$ belongs to a broader class of ``asymptotically $\psi$-locally
linear functions'', i.e.\ such functions that, for any fixed  $v$ and $x\to\infty$,
$$
F_+ (t+v\psi(t) )=F_+(t) (1-cv+o(1) ),\quad c={\rm const}>0.
$$
However, such a representation is impossible as $1-cv<0$ when $v>1/c$.

\section{The characterization of $\psi$-l.c.f.'s}

The aim of the present section is to prove that, for any $\psi$-l.c.f.\ $g$, convergence
\eqref{ap5}  is uniform in~$v$ on any compact set and, moreover, that $g$ admits an
integral representation similar to~\eqref{ap2}, \eqref{ap3'}. To do that, we will need
some restrictions on the function~$\psi$.

We assume that $\psi$ is a non-decreasing function such that $\psi(x)=o(x)$ as
$x\to\infty$. For such functions, we introduce the following condition:

\medskip
{\bf (A)} \ {\em For any fixed $v>0$, there exists a value $a(v)\in(0,\infty)$ such that
\begin{equation}
 \label{ap*_}
  \frac{\psi (x-v\psi(x) )}{\psi(x)}\geq
a(v) \quad \mbox{for all large enough \ $x$}.
\end{equation}
}
Of course, one can assume without loss of generality that $a(v)$ is non-increasing.

Now note that, letting $y:= x+v  \psi (x)>x  $ and   using the monotonicity of $\psi$,
one has
\[
\psi (y- v  \psi (y)) \le \psi (x).
\]
Therefore, relation~\eqref{ap*_} implies that, for all large enough~$x$,
\[
  \frac{\psi (x + v \psi (x))}{\psi (x)}\equiv \frac {\psi (y)}{\psi (x)}
   \le  \frac{\psi (y)}{\psi (y- v  \psi (y))}
   \le \frac1{a(v)}\in (0,\infty).
\]
Thus, any function  $\psi$ satisfying condition {\bf (A)} will aslo satifsy the
following relation: for any fixed $v>0$,
\begin{equation}
 \label{ap*_b}
\frac{\psi (x+v\psi(x) )}{\psi(x)}\leq \frac{1}{a(v)} \quad \mbox{for all large enough \
$x$}.
\end{equation}
Observe that the converse is not true: it is not hard to construct an example of a
(piece-wise linear, globally Lipschitz) non-decreasing function  $\psi$ which satisfies
a condition of the form~\eqref{ap*_b}, but for which condition {\bf (A)} will hold for
no finite function~$a(v)$.

It is clear that if $\psi$ is a $\psi$-l.c.f., $\psi(x)=o(x)$, then $\psi$ satisfies
condition~{\bf (A)}.

Introduce class $\mathcal{K}$ consisting of non-decreasing  functions $\psi(x)\ge 1$,
$x\ge 0,$ that satisfy   condition~{\bf (A)} for a function $a(v)$  such that
\begin{equation}
 \label{14}
\int_0^\infty  a(u) \, du = \infty.
\end{equation}

Class $\mathcal{K}_1$ we define as the class of continuous  r.v.f.'s of index $\alpha<1$
and such that  $x/\psi(x)\uparrow\infty$ as $x\to\infty$ and the following ``asymptotic
smoothness'' condition is met:
\begin{equation}
 \label{ap9}
\psi(x+\Delta)=\psi(x)+\frac{\alpha\Delta\psi(x)}{x}\, (1+o(1) ) \quad\mbox{as }
x\to\infty,\quad\Delta=v\psi(x),\quad v={\rm const}.
\end{equation}
Clearly, $\mathcal{K}_1\subset\mathcal{K}$. Condition~\eqref{ap9} is always satisfied
for any $\Delta\geq c_1={\rm const}$, $\Delta=o(x)$, provided that the function $L(x)$
is differentiable and $L'(x)=o\left( {L(x)}/{x}\right)$ as $x\to\infty$.

In the assertions to follow, it will be assumed that $\psi$ belongs to $\mathcal{K}$ or
$\mathcal{K}_1$. We will not dwell on how far  the conditions $\psi\in\mathcal{K}$ or
$\psi\in\mathcal{K}_1$ can be extended. The function $\psi$ specifies the ``asymptotic
local constancy zone width'' of the function~$g$ under consideration, and what matters
for us is just the growth rate of~$\psi(x)$ as $x\to\infty$. All the other properties of
$\psi$ (its smoothness, the presence of oscillations etc.)  are for us to choose, and so
we can assume the function $\psi$ to be as  ``smooth'' as we need.  In this sense, the
assumption that $\psi$ belongs to the class  $\mathcal{K}$ or~$\mathcal{K}_1$ is not
restrictive. For example, it is quite natural to assume in Theorems~B and~C from Section
2 that~$\psi\in\mathcal{K}_1$.

The following assertion extends the property {\bf($\boldsymbol{\rm U_1}$)} of l.c.f.'s
to $\psi$-l.c.f.'s.

\begin{Theorem}
 \label{ath0}
If\, $g$ is a $\psi$-l.c.f.\ with $\psi\in\mathcal{K}$, then convergence  in~\eqref{ap5}
is uniform: for any fixed real numbers $v_1<v_2$,
\begin{equation}
\label{ap5_1}
{\bf (U_{\boldsymbol{\psi}})}\hspace{32mm}
 \lim_{x\to \infty} \sup_{v_1\le v \le v_2}
  \bigg| \frac{g \big(x+v\psi(x)\big)}{g(x)}-1\bigg|=0.
 \hspace{30mm}
\end{equation}
\end{Theorem}

Observe that, for monotone~$g$, the condition $\psi\in\mathcal{K}$ in Theorem~\ref{ath0}
is superfluous. Indeed, assume for definiteness that $g$ is a non-decreasing
$\psi$-l.c.f. Then, for any   $v$ and $v(x)\to v$, there is a   $v_0>v$ such that, for
all sufficiently large~$x$, one has $v(x)<v_0$, and therefore
\begin{equation}
\label{ap11'}
\limsup_{x\to\infty} \,\frac{g (x+v(x)\psi(x) )}{g(x)}
  \leq
\limsup_{x\to\infty} \,\frac{g (x+v_0\psi(x) )}{g(x)}=1.
\end{equation}
A converse inequality for $\liminf$ is established in a similar way. As a consequence,
\begin{equation}\label{ap12'}
\lim_{x\to\infty}\frac{g (x+v(x)\psi(x) )}{g(x)}=1,
\end{equation}
which is easily seen to be equivalent to~\eqref{ap5_1} (cf.~\eqref{ap4}).

It is not hard to see that monotonicity property required to derive~\eqref{ap11'},
\eqref{ap12'}, could be somewhat relaxed.

Now set
\begin{equation}
\label{ap12'_}
\gamma(x):=\int_1^x\frac{dt}{\psi(t)}.
\end{equation}

\begin{Theorem}
 \label{ath1}
Let  $\psi\in\mathcal{K}$. Then $g$ is a $\psi$-l.c.f.\ iff it admits a representation
of the form
\begin{equation}
\label{ap10}
\boldsymbol{\rm (I_\psi)}
 \hspace{30mm}
 g(x)= c(x) \exp\bigg\{\int_1^{e^{\gamma(x)}}\frac{\varepsilon(t)}{t}\,dt\bigg\},\qquad
x\geq 1,
 \hspace{26mm}
\end{equation}
where $c(t)$ and $\varepsilon(t)$ have the same properties as in~{\bf (I)}.
\end{Theorem}

Since, for any  $\varepsilon>0$ and all large enough~$x$,
$$
\int_1^{e^{\gamma(x)}}\frac{\varepsilon(t)}{t}\,dt
 < \varepsilon\ln e^{\gamma(x)}=\varepsilon\gamma(x)
$$
and a similar lower bound holds true, Theorem~\ref{ath1} implies the following result.

\begin{Corollary}
 \label{aco1}
If\/ $\psi\in\mathcal{K}$ and $g$ is a $\psi$-l.c.f., then
$$
g(x)=e^{o(\gamma(x))},\qquad x\to\infty.
$$
\end{Corollary}

For $\psi\in\mathcal{K}_1$ we put
$$
\theta(x):=\frac{x}{\psi(x)}.
$$
Clearly,  $\theta(x)\sim(1-\alpha)\gamma(x)$ as $x\to\infty$.

\begin{Theorem}
 \label{ath2}
Let  $\psi\in\mathcal{K}_1$. Then the assertion of Theorem~$\ref{ath1}$ holds true with
$\gamma(x)$ replaced by~$\theta(x)$.
\end{Theorem}

\begin{Corollary}
 \label{aco2}
If $\psi\in\mathcal{K}_1$ and $g$ is a $\psi$-l.c.f., then
$$
g(x)=e^{o(\theta(x))},\qquad x\to\infty.
$$
\end{Corollary}

Since the function $\theta(x)$ has a ``more explicit'' representation in terms of~$\psi$
than the function~$\gamma(x)$, the assertions of Theorem~\ref{ath2} and
Corollary~\ref{aco2}  display the asymptotic properties $\psi$-l.c.f.'s in a more
graphical way than those of Theorem~\ref{ath1} and Corollary~\ref{aco1}. A deficiency of
Theorem~\ref{ath2} is the fact that the condition $\psi\in\mathcal{K}_1$ is more
restrictive than the condition that~$\psi\in\mathcal{K}$. It is particularly essential
that, in the former condition, the equality $\alpha=1$ is excluded for the
index~$\alpha$ of the r.v.f.~$\psi$.

\section{Proofs}

{\em Proof of Theorem~\ref{ath0}}. Our proof will use an argument modifying H.~Delange's
proof of property~{\bf (U)} (see e.g.\ p.\,6 in~\cite{bingham}) or \S~1.1
in~\cite{bor_bor}).

Let $l(x): =\ln g(x).$ It is clear that~\eqref{ap5} is equivalent to the convergence
\begin{equation}
 \label{ap5_+}
l(x+v\psi(x)) - l(x) \to 0, \qquad x\to \infty,
\end{equation}
for any fixed~$v\in \mathbb{R}$. To prove the theorem, it suffices to show that
$$
H_{v_1, v_2} (x) := \sup_{v_1 \le v\le v_2} \bigl|l(x+v\psi(x)) -
l(x)\bigr|  \to 0, \qquad x\to \infty.
$$
It is not hard to see that the above relation will follow from the convergence
\begin{equation}
 \label{ap5_++}
H_{0, 1} (x)  \to 0, \qquad x\to \infty.
\end{equation}
Indeed, let   $v_1<0$ (for  $v_1\ge 0$ the argument will be even simpler) and
\[
x_0:= x+ v_1\psi (x),\qquad x_k := x_0 + k \psi (x_0), \qquad
k=1,2,\ldots
\]
By virtue of condition~{\bf (A)}, one has  $\psi (x_0) \ge a(-v_1) \psi (x)$ with
$a(-v_1)>0$. Therefore, letting $n:= \lfloor (v_2-v_1)/ a(-v_1)\rfloor +1,$  where
$\lfloor x\rfloor$ denotes the integer part of~$x$, we obtain
\[
H_{v_1, v_2} (x) \le \sum_{k=0}^n H_{0, 1} (x_k),
\]
which establishes the required implication.

Assume without loss of generality that $\psi (0)=1.$ To prove~\eqref{ap5_++}, fix an
arbitrary small $\varepsilon\in \bigl(0, a(1)/(1+ a(1)) \bigr)$ and set
\begin{align*}
I_x      := [x, x & + 2\psi (x)],   \qquad
 I^*_x   := \{ y\in I_x : \, |l(y) - l(x) | \ge \varepsilon /2\} ,\\
 & I^*_{0,x}  := \{ u\in I_0 : \, |l(x+u \psi (x)) - l(x)|  \ge
\varepsilon /2\}.
\end{align*}
One can easily see that all these sets are measurable and
\[
 I^*_x  = x + \psi (x) I^*_{0,x},
\]
so that for the Lebesgue measure $\mu (\cdot )$ on $\mathbb{R}$  we have
\begin{equation}\label{ap5_+++}
\mu(I^*_x)  =   \psi (x) \mu (I^*_{0,x}).
\end{equation}
It follows from~\eqref{ap5_+} that, for any  $u\in I_0$, the value of the indicator
${\bf 1}_{I^*_{0,x}} (u)$ tends to zero as~$x\to\infty$. So, by the dominated
convergence theorem,
\[
\int_{I_0} {\bf 1}_{I^*_{0,x}} (u) \, du \to 0, \qquad x\to 0.
\]
From here and \eqref{ap5_+++} we see that there exists an $x_{(\varepsilon)}$ such that
\[
\mu(I^*_x)  \le \frac{\varepsilon}{2}\, \psi (x), \qquad x\ge
x_{(\varepsilon)}.
\]
Now observe that, for any $s\in [0,1]$, the set $I_x \cap I_{x+ s\psi (x)} = [x+ s\psi
(x), x + 2 \psi (x)]$ has length  $(2-s) \psi (x) \ge \psi (x).$ Hence,  for $x\ge
x_{(\varepsilon)}$, the set
\[
J_{x, s} := \bigl( I_x \cap I_{x+ s\psi (x)}\bigr) \setminus \bigl(
I^*_x \cup I^*_{x+ s\psi (x)}\bigr)
\]
will have the length
\begin{align*}
\mu (J_{x, s} )
 &  \ge \psi (x) - \frac{\varepsilon}2 \, \bigl[\psi (x) + \psi (  x + s\psi (x) )\bigr] \\
 &  \ge \psi (x) - \frac{\varepsilon}2 \, \biggl(1 + \frac1{a(1)}\biggr)
\psi (x) \ge \frac12 \psi (x) \ge \frac12,
\end{align*}
where we used relation~\eqref{ap*_b} to establish the second  inequality. Therefore
$J_{x,s}\neq \varnothing  $ and one can choose a point $y\in J_{x,s}$. Then $y \not \in
I^*_x$ and $y \not \in I^*_{x+ s\psi (x)}$, so that
\[
|l (x+ s\psi (x)) - l (x) | \le |l (x+ s\psi (x)) - l (y) | + |
l(y)- l (x) | <\varepsilon.
\]
Since this relation holds for any $s\in [0,1]$, the required convergence \eqref{ap5_++}
and hence the assertion of Theorem~\ref{ath0} are proved.\hfill$\square$

\medskip

{\em Proof of Theorem~\ref{ath1}.} First let $g$ be a $\psi$-l.c.f.\ with
$\psi\in\mathcal{K}$. Since $\psi(t)=o(t)$, one has $\gamma(x)\uparrow \infty$ as
$x\uparrow\infty$ (see~\eqref{ap12'_}). Moreover, the function  $\gamma(x)$ is
continuous and so always has an inverse  $\gamma^{(-1)}(t)\uparrow \infty$ as
$t\to\infty$, so that we can consider the composition function
\[
g_\gamma (t):=(g\circ \gamma^{(-1)}) (t).
\]
If we show that $g_\gamma$ is an l.c.f.\ then representation~\eqref{ap10} will
immediately follow from the relation  $  g(x) = g_\gamma (\gamma(x) )$ and property~{\bf
(I$_{\bf 1}$)}.

By virtue of the uniformity property ${\bf (U_{\boldsymbol \psi})}$  which holds for $g$
by Theorem~\ref{ath0}, for any bounded  function  $r(x)$ one has
\begin{equation}
 \label{ap1.9.1}
g_\gamma\big(\gamma(x)\big)
 \equiv
    g(x)\sim g\big(x+r(x)\psi(x)\big)
 =
    g_\gamma\big(\gamma (x + r(x)\psi(x) )\big).
\end{equation}

Next we will show that, for a given $v$ (let $v>0$ for definiteness), there is a bounded
(as $x\to\infty$) value $r(x,v)$ such that
\begin{equation}
 \label{ap2.9.1}
\gamma (x+r(x,v)\psi(x) )=\gamma(x)+v.
\end{equation}
Indeed, we have
$$
\gamma (x+r\psi(x) )-\gamma(x)= \int_x^{x+r\psi(x)}\frac{dt}{\psi(t)} =
\int_0^r\frac{\psi(x)\, dz}{\psi(x+z\psi(x))} =:I(r,x),
$$
where, by Fatou's lemma and relation~\eqref{ap*_b},
$$
 \liminf_{x\to\infty}I(r,x)\geq\int_0^r\liminf_{x\to\infty}
\frac{\psi(x)}{\psi(x+z\psi(x))}\, dz\geq I(r):= \int_0^r a(z)dz\uparrow\infty
$$
as $r\uparrow\infty$ (see~\eqref{14}). Since, moreover, for any~$x$ the function
$I(r,x)$ is continuous in~$r$, there exists an $r(v,x)\leq r_v<\infty$ such that
$I\big(r(v,x),x\big)=v$, where $r_v$ is the solution of the equation~$I(r)=v$.

Now choosing the $r(x)$ in~\eqref{ap1.9.1} to be the function $r(x,v)$
from~\eqref{ap2.9.1} we obtain that
$$
g_\gamma\big(\gamma(x)\big)\sim g_\gamma\big(\gamma(x)+v\big)
$$
as $x\to\infty$, which means that  $g_\gamma$  is an l.c.f.\ and hence~\eqref{ap10}
holds true.

Conversely, let representation~\eqref{ap10} hold true. Then, for a fixed~$v\geq 0$, any
$\varepsilon>0$ and $x\to\infty$, one has
\begin{align}
 \left|\ln\,\frac{g\big(x+v\psi(x)\big)}{g(x)}\right|
  & \le
 \int_{e^{\gamma(x)}}^{e^{\gamma(x+v\psi(x))}}
 \frac{\big|\varepsilon(t)\big|}{t}\,dt+o(1)
  \le
 \Big(\gamma\big(x+v\psi(x)\big)-\gamma(x)\Big)\varepsilon+o(1)
  \nonumber\\
 &\leq
 \varepsilon\int_0^v\frac{\psi(x)ds}{\psi(x+s\psi(x))}+o(1)
  \leq\varepsilon v+o(1).
 \label{ap1.9.2}
\end{align}
This clearly means that the left-hand side of this relation is  $o(1)$ as $x\to\infty$.

If $v=-u<0$ then, bounding in a similar fashion the integral
$$
\int_{e^{\gamma(x-u\psi(x))}}^{e^{\gamma(x)}}\frac{\big|\varepsilon(t)\big|dt}{t}
 \leq \varepsilon\int_0^u\frac{\psi(x)ds}{\psi(x-s\psi(x))},
$$
we will obtain from condition {\bf (A)} that
$$
\limsup_{x\to\infty}\left|\ln\frac{g(x+v(\psi(x)))}{g(x)}\right|
 \leq \varepsilon\int_0^u \limsup_{x\to\infty}\frac{\psi(x)ds}{\psi(x-s\psi(x))}
 \leq \varepsilon\int_0^u\frac{ds}{a(s)},
$$
so that the left-hand side of \eqref{ap1.9.2} is still $o(1)$ as $x\to\infty$. Therefore
$g (x+v\psi(x) )\sim g(x)$ and hence $g$ is a $\psi$-l.c.f. Theorem~\ref{ath1} is
proved. \hfill$\square$

\medskip

It is evident that the assertion of Theorem~\ref{ath1} can also be stated as follows:
for $\psi\in\mathcal{K}$, a function  $g$ is a $\psi$-l.c.f.\ iff  $g_\gamma(x)$ is an
l.c.f.\ (which, in turn, holds iff $g_\gamma(\ln x)$ is an s.v.f.).

\medskip
{\em Proof of Theorem~\ref{ath2}}. One can employ an argument similar to the one used to
prove Theorem~\ref{ath1}.

Since the function $\theta(x)$ is continuous and increasing, it has an
inverse~$\theta^{(-1)}(t)$. It is not hard to see that if  $\psi$ обладает has
property~\eqref{ap9}, then the function  $\theta(x)=x/\psi(x)$ also possesses a similar
property: for a fixed $v$ and $\Delta=v\psi(x)$, $x\to\infty$, one has
\[
\theta(x+\Delta)=\theta(x) + \frac{(1-\alpha)\Delta\theta(x)}{x}\,  (1+o(1) ).
\]
Therefore, as $x\to\infty$,
$$
\theta (x+v\psi(x) )=\theta(x)+(1-\alpha)v (1+o(1)) .
$$
As the function  $\theta$ is monotone and continuous, this relation means that, for
any~$v$, there is a function $v(x) \to v$ as $x\to\infty$ such that
\begin{equation}
 \label{ap16}
\theta (x+v(x)\psi(x) )=\theta(x)+(1-\alpha)v.
\end{equation}

Let  $g$ be a $\psi$-l.c.f. Then, for the function $g_\theta:=g\circ \theta^{(-1)},$ we
obtain by virtue of~\eqref{ap16} that
$$
g_\theta (\theta(x) )\equiv g(x)\sim g (x+v(x)\psi(x) ) =
g_\theta\left(\theta (x+v(x)\psi(x) )\right)= g_\theta
(\theta(x)+(1-\alpha)v ).
$$
Since $\theta(x)\to\infty$ as $x\to\infty$, the above means that  $g_\theta$ is an
l.c.f. The direct assertion of the integral representation theorem follows from here
and~\eqref{ap3'}.

The converse assertion is proved in the same way as in  Theorem~\ref{ath1}.
Theorem~\ref{ath2} is proved. \hfill$\square$

\medskip

Similarly to our earlier argument, it follows from Theorem~\ref{ath2} that if
$\psi\in\mathcal{K}_1$ then  $g$ is a $\psi$-l.c.f.\ iff $g_\theta$ is an l.c.f.\ (and
$g_\theta(\ln x)$ is an s.v.f.).

\medskip

{\bf Acknowledgements.} Research supported by the Russian Federation President    Grant
NSh-3695.2008.1, Russian Foundation for Basic Research Grant 08--01--00962  and the ARC
Centre of Excellence for Mathematics and Statistics of Complex Systems.


\begin{thebibliography}{X}

\bibitem{bingham}
{\sc Bingham, N.\,H., Goldie, C.\,M.  and Teugels, I.\,L.} (1987) {\em Regular
variation}. Cambridge: Cambridge University Press.

\bibitem{bor_bor}
{\sc  Borovkov, A.\,A.  and Borovkov, К.\,А.} (2008) {\em Asymptotic analysis of random
walks: Heavy-tailed distributions.} Cambridge: Cambridge University Press.

\bibitem{Cline-Hsing} {\sc Cline, D.\,B.\,H. and Hsing, T.} (1991) Large deviations probabilities for
sums and maxima of random variables with heavy or subexponential tails. Preprint: Texas
A\,\&M University.

\bibitem{Cline}
{\sc Cline, D.\,B.\,H.} Intermediate regular and $\Pi$ variation (1994) {\em Proc.
London Math. Soc.} {\bf 68,} 594--616.

\bibitem{Ng}
{\sc Ng, K.\,W., Tang, Q., Yan, J.-A. and Yang H.} (1994) Precise large deviation for
sums of random variables with consistently varying tails. {\em J. Appl. Prob. } {\bf
41,} 93--107.

\end{thebibliography}
\end{document}